\documentclass{amsart}
\usepackage{amsmath, amsthm, amsfonts}
\usepackage{verbatim}

\numberwithin{equation}{section}
\newtheorem{lemma}{Lemma}[section]

\newtheorem{theorem}{Theorem}[section]

\newcommand{\lgn}{\frac{\log{\log{n}}}{n}}
\newcommand{\elg}{\frac{1}{\log{n}}}
\newcommand{\Exp}{\textup{E}}

\newcommand{\bfl}{\begin{flalign*}}
\newcommand{\efl}{\end{flalign*}}
\newcommand{\cont}{C\left([-\pi,\pi]\right)}
\newcommand{\diff}{C^1\left([-\pi,\pi]\right)}

\begin{document}

\title[Real zeros of a random algebraic polynomial]{The real zeros of a random algebraic polynomial with dependent coefficients}
\author[J. Matayoshi]{Jeffrey Matayoshi}
\address{Department of Mathematics\\
340 Rowland Hall\\
University of California, Irvine\\
Irvine, CA 92697-3875}
\email{jsmatayoshi@gmail.com}
\thanks{This research was partially supported by NSF grant DMS-0706198}

\begin{abstract}
Mark Kac gave one of the first results analyzing random polynomial zeros.  He considered the case of independent standard normal coefficients and was able to show that the expected number of real zeros for a degree $n$ polynomial is on the order of $\frac{2}{\pi}\log{n}$, as $n\rightarrow \infty$.  Several years later, Sambandham considered two cases with some dependence assumed among the coefficients.  The first case looked at coefficients with an exponentially decaying covariance function, while the second assumed a constant covariance.  He showed that the expectation of the number of real zeros for an exponentially decaying covariance matches the independent case, while having a constant covariance reduces the expected number of zeros in half.  In this paper we will apply techniques similar to Sambandham's and extend his results to a wider class of covariance functions.  Under certain restrictions on the spectral density, we will show that the order of the expected number of real zeros remains the same as in the independent case.
\end{abstract}

\maketitle

\section{Introduction}

One of the earliest results on the expected number of real zeros of the random polynomial given by
\begin{equation} \label{poly}
P_n(x)=\sum_{k=0}^{n}{X_{k}x^k}
\end{equation}
came from Mark Kac \cite{Kac}.  Kac considered the case when the coefficients are assumed to be independent standard normal random variables and was able to show that the value of the expected number of zeros is on the order of $\frac{2}{\pi}\log{n}$, as $n\rightarrow \infty$.  More recently, Edelman and Kostlan \cite{EdelKost} derived a similar result, but in doing so they gave a nice geometric argument and derived formulas that hold for a wider class of coefficients.  

A natural generalization of this problem is to assume some dependence among the coefficients.  Let $X_0, X_1, \ldots$ be a stationary  sequence of normal random variables, where the covariance function is given by
\[
\Gamma(k)=\textup{E}[X_0X_k], \quad \Gamma(0)=1.
\]
Under these assumptions, two important results came from Sambandham.  The first assumes that $\Gamma(k)=\rho^{k}$, where $\rho\in(0,\frac{1}{2})$ \cite{Sam77}.  In this case, it was shown that the expected number of zeros is on the same order as when the coefficients are independent.  The second result assumes the covariance function is constant; that is, $\Gamma(k)=\rho$ for any $k$, where $\rho \in (0,1)$ \cite{Sam76,Sam79}.  Here, it was shown that the order of the expected number of zeros is only $\frac{1}{\pi}\log{n}$ as $n\rightarrow \infty$, or half as many as before.  The intuitive reasoning for this is as follows.  For the first result, the exponential decay of the covariance function is fast enough to ``simulate'' independent behavior.  In other words, the dependence between the coefficients disappears at a rate fast enough to have no effect on the behavior of the zeros.   For the second result, however, the constant covariance results in the coefficients tending to have the same sign.  This causes most of the positive zeros to disappear, leaving only half as many zeros as before.  

A further question that could be asked is whether or not this result extends to a wider class of covariance functions, where the decay rates are between those considered by Sambandham.  In what follows, we will show that the same asymptotic value holds for the expected number of real zeros when the spectral density (which will be introduced below) is positive and continuous.  If we assume further that the spectral density also has one continuous derivative, we will be able to derive an explicit value for the order of the error term.  Noting that absolute summability of the covariance function guarantees the continuity of the spectral density, we will then see that behavior similar to the independent case can be expected for covariance functions with a wide range of decay rates.  The main result is stated as follows.

\begin{theorem} \label{roots}
Let $P_n(x)$ be the polynomial given in \eqref{poly}, where the coefficients $X_0, X_1, \ldots$ form a stationary sequence of standard normals, with covariance function $\Gamma(k)$ and spectral density $f(\phi)$.  Assume that the spectral density does not vanish. Letting $N(\alpha,\beta)$ be the number of zeros of $P_n(x)$ in the interval $(\alpha,\beta)$, it follows that
\begin{flalign*}
\textup{E}[N(-\infty, \infty)] &\sim \frac{2}{\pi}\log{n}, \qquad \mbox{for $f\in\cont$},\\
\textup{E}[N(-\infty, \infty)] &=    \frac{2}{\pi}\log{n} + O\left(\log{\log{n}}\right), \qquad \mbox{for $f\in\diff$},\\
\end{flalign*}
as $n\rightarrow \infty$.
\end{theorem}

Before we proceed any further, a comment must be made.  As mentioned on page 83 of \cite{Sam77}, if we consider the function
\begin{flalign*}
x^nP_n\left(\frac{1}{x}\right) &= x^n\left(X_0 + X_1\frac{1}{x} + \cdots + X_n\frac{1}{x^n}\right)\\
										&= X_0x^n + X_1x^{n-1} + \cdots + X_n,
\end{flalign*}			
it can be seen that whenever there is a zero of $x^nP_n\left(\frac{1}{x}\right)$ in $(1,\infty)$, there is also a zero of $P_n(x)$ in $(0,1)$.  Thus, since the distribution of the zeros of the two functions are the same, it is sufficient to only look at the interval from $(0,1)$.  A similar argument works for the negative real line and allows us to restrict our analysis to the interval $(-1,1)$.  By then taking twice the result, we will have a value for the total expected number of real zeros.
 
Our proof will follow that of Sambandham, with some necessary modifications to account for the more general assumptions made on the coefficients.  The first step is to show that there is a negligible amount of zeros on the intervals $(0,1-\elg)$, $(1-\frac{\log{\log{n}}}{n},1)$, $(-1+\elg,0)$ and $(-1,-1+\frac{\log{\log{n}}}{n})$.  Following this, we will then show that the number of real zeros in the intervals $(1-\elg,1-\frac{\log{\log{n}}}{n})$ and $(-1+\frac{\log{\log{n}}}{n},-1+\elg)$ are each on the order of $\frac{1}{2\pi}\log{n}$.  One significant difference from Sambandham's work is that rather than approximating the number of zeros with a specially chosen function, we will instead derive the asymptotic behavior directly from the Kac-Rice formula using the spectral density formulation of the covariance function (see below).  The result is that the number of zeros in $(-1,1)$ is on the order of $\frac{1}{\pi}\log{n}$ and, from the comments above, it follows that the total expected number of real zeros is on the order of $\frac{2}{\pi}\log{n}$.

\section{Preliminary Results}

We will now discuss some results that are variously attributed to Bochner, Herglotz, and Khinchine (see \cite{Brei,CramLead}).  Any covariance function, $\Gamma(k)$, can be expressed as
\[
\textup{E}[X_0X_k]=\Gamma(k)=\int_{-\pi}^{\pi}e^{-ik\phi}F(d\phi),
\]
where $F(\phi)$ is real, never-decreasing, and bounded.  Furthermore, if $F(\phi)$ is also absolutely continuous, we have the formula
\begin{equation} \label{cov}
\Gamma(k)=\int_{-\pi}^{\pi}e^{-ik\phi}f(\phi)d\phi,
\end{equation}
where $f(\phi)$ is called the spectral density of the covariance function.  A sufficient condition for the existence of $f(\phi)$ is that $\Gamma(k)$ is absolutely summable.  Additionally, in this case it will be non-negative, continuous, and of the form
\[
f(\phi) = \frac{1}{2\pi}\sum_{k=-\infty}^{\infty}\Gamma(k)e^{ik\phi}.
\]

Now, recalling that $N(\alpha,\beta)$ is the number of zeros of $P_n(x)$ in the interval $(\alpha,\beta)$, the Kac-Rice formula \cite{Farahmand} gives the expected number of zeros on this interval as 
\begin{equation} \label{Kac}
\textup{E}[N(\alpha,\beta)]=\frac{1}{\pi}\int_{\alpha}^{\beta}{\frac{\sqrt{AC-B}}{A}dx},
\end{equation}
where 
\begin{equation} \label{ABC}
\begin{split}
A(x) &= \textup{E}\left[P_n^2(x)\right] = \sum_{k=0}^{n}{\sum_{j=0}^{n}{\Gamma(k-j)x^{k+j}}},\\
B(x) &= \textup{E}\left[P_n(x)P_n'(x)\right] = \sum_{k=0}^{n}{\sum_{j=0}^{n}{\Gamma(k-j)kx^{k+j-1}}},\\
C(x) &= \textup{E}\left[(P_n'(x))^2\right] = \sum_{k=0}^{n}{\sum_{j=0}^{n}{\Gamma(k-j)kjx^{k+j-2}}}.\\
\end{split}
\end{equation}
Applying \eqref{cov} we can rewrite these as
\begin{equation} \label{ABC2}
\begin{split}
A &= \int_{-\pi}^{\pi}{\sum_{k=0}^{n}{\sum_{j=0}^{n}{e^{-i(k-j)\phi}x^{k+j}f(\phi)d\phi}}},\\
B &= \int_{-\pi}^{\pi}{\sum_{k=0}^{n}{\sum_{j=0}^{n}{e^{-i(k-j)\phi}kx^{k+j-1}f(\phi)d\phi}}},\\
C &= \int_{-\pi}^{\pi}{\sum_{k=0}^{n}{\sum_{j=0}^{n}{e^{-i(k-j)\phi}kjx^{k+j-2}f(\phi)d\phi}}}.\\
\end{split}
\end{equation}
We are now ready to prove our first lemma.
\begin{lemma} \label{bound}
For the intervals $(-1,-1+\frac{\log{\log{n}}}{n})$, $(-1+\elg,0)$, $(0,1-\elg)$, and $(1-\frac{\log{\log{n}}}{n},1)$, the expected number of zeros is $O(\log{\log{n}})$.
\end{lemma}
\begin{proof}
Following the method of Sambandham, we will define the function
\begin{equation} \label{H}
\begin{split}
H(x,y) &= \sum^{n}_{k=0}{\sum^{n}_{j=0}{e^{-i(k-j)\phi}x^{k}y^{j}}}\\
       &= \sum_{k=0}^{n}{e^{-ik\phi}x^{k}\sum_{j=0}^{n}{e^{ij\phi}y^{j}}}\\
       &= \frac{1-x^{n+1}e^{-i(n+1)\phi}}{1-xe^{-i\phi}} \cdot \frac{1-y^{n+1}e^{i(n+1)\phi}}{1-ye^{i\phi}},
\end{split}
\end{equation}
to assist in our computations.  As done on page 85 of \cite{Sam77}, plugging into our formula for $A$ leads to the expression 
\begin{equation} \label{intA}
\begin{split}
A &= \int_{-\pi}^{\pi}H(x,x)f(\phi)d\phi\\
  &= \int_{-\pi}^{\pi}\frac{1-x^{n+1}e^{-i(n+1)\phi}}{1-xe^{-i\phi}} \cdot \frac{1-x^{n+1}e^{i(n+1)\phi}}{1-xe^{i\phi}}f(\phi)d\phi.\\
\end{split}
\end{equation}
Again following \cite{Sam77}, we also get
\begin{equation} \label{intB}
\begin{split}
B &= \int_{-\pi}^{\pi}\left[\frac{\partial H(x,y)}{\partial y}\right]_{y=x}f(\phi)d\phi\\
  &= \int_{-\pi}^{\pi}\left(\frac{1-x^{n+1}e^{-i(n+1)\phi}}{1-xe^{-i\phi}}\right)\\
  &  \quad \cdot \left(\frac{-(n+1)x^ne^{i(n+1)\phi}\left(1 - xe^{i\phi}\right) - \left(1 - x^{n+1}e^{i(n+1)\phi}\right)\left(-e^{i\phi}\right)}
  {\left(1-xe^{i\phi}\right)^2}\right)f(\phi)d\phi,\\
\end{split}
\end{equation}
and
\begin{equation} \label{intC}
\begin{split}
C &= \int_{-\pi}^{\pi}\left[\frac{\partial^2H(x,y)}{\partial x\partial y}\right]_{y=x}f(\phi)d\phi\\
  &= \int_{-\pi}^{\pi}\left(\frac{-(n+1)x^ne^{-i(n+1)\phi}\left(1 - xe^{-i\phi}\right) - \left(1 - x^{n+1}e^{-i(n+1)\phi}\right)\left(-e^{-i\phi}\right)}{\left(1-xe^{-i\phi}\right)^2}\right)\\
  &  \quad \cdot \left(\frac{-(n+1)x^ne^{i(n+1)\phi}\left(1 - xe^{i\phi}\right) - \left(1-x^{n+1}e^{i(n+1)\phi}\right)\left(-e^{i\phi}\right)}
  {\left(1-xe^{i\phi}\right)^2}\right)f(\phi)d\phi.\\
\end{split}
\end{equation}
For $x\in(0, 1-\elg)$ we have
\begin{flalign*}
A &\sim \int_{-\pi}^{\pi}\frac{1}{(1-xe^{-i\phi})(1-xe^{i\phi})}f(\phi)d\phi,\\
\end{flalign*}
and
\begin{flalign*}
C &\sim \int_{-\pi}^{\pi}\frac{1}{(1-xe^{-i\phi})^2(1-xe^{i\phi})^2}f(\phi)d\phi\\
  &\leq    \frac{1}{(1-x)^2}\int_{-\pi}^{\pi}\frac{1}{(1-xe^{-i\phi})(1-xe^{i\phi})}f(\phi)d\phi.\\
\end{flalign*}
If we consider the quotient $\frac{\sqrt{AC-B^2}}{A}$, then
\[
\frac{\sqrt{AC-B^2}}{A} < \left(\frac{C}{A}\right)^{1/2} \leq \frac{c}{1-x}.
\]
Plugging into \eqref{Kac}, it follows that
\begin{equation} \label{int1}
\textup{E}\left[N\left(0,1-\elg\right)\right] = O\left(\log{\log{n}}\right).
\end{equation}

To handle the interval $(-1+\elg,0)$ we will substitute in $-x$, where \linebreak[4]$x\in(0, 1-\elg)$, to get
\begin{flalign*}
A &=       \int_{-\pi}^{\pi}H(-x,-x)f(\phi)d\phi\\
  &=       \int_{-\pi}^{\pi}\frac{1-(-x)^{n+1}e^{-i(n+1)\phi}}{1+xe^{-i\phi}} \cdot \frac{1-(-x)^{n+1}e^{i(n+1)\phi}}{1+xe^{i\phi}}f(\phi)d\phi\\
  &\sim \int_{-\pi}^{\pi}\frac{1}{(1+xe^{-i\phi})(1+xe^{i\phi})}f(\phi)d\phi.\\
\end{flalign*}
Likewise,
\begin{flalign*}
C &\sim \int_{-\pi}^{\pi}\frac{1}{\left(1+xe^{-i\phi}\right)^2\left(1+xe^{i\phi}\right)^2}f(\phi)d\phi\\
  &\leq \frac{1}{(1-x)^2}\int_{-\pi}^{\pi}\frac{1}{\left(1+xe^{-i\phi}\right)\left(1+xe^{i\phi}\right)}f(\phi)d\phi.
\end{flalign*}
Once more we have the inequality
\[
\frac{\sqrt{AC-B^2}}{A} < \left(\frac{C}{A}\right)^{1/2} \leq \frac{c}{1-x}.
\]
Applying \eqref{Kac}, we then have
\begin{equation} \label{int2}
\textup{E}\left[N\left(-1+\elg,0\right)\right] = O\left(\log{\log{n}}\right).
\end{equation}

Finally, we will consider the intervals $(-1, -1+\frac{\log{\log{n}}}{n})$ and $(1-\frac{\log{\log{n}}}{n}, 1)$.  Similar to inequality derived on page 87 in \cite{Sam77}, for $x\in(1-\frac{\log{\log{n}}}{n}, 1)$ we have
\begin{equation*}
\frac{\sqrt{AC-B^2}}{A} < \left(\frac{C}{A}\right)^{1/2} < \left(\frac{n^2\sum^{n}_{k=0}{\sum^{n}_{j=0}{\Gamma(k-j)x^{k+j-2}}}}
		{x^2\sum^{n}_{k=0}{\sum^{n}_{j=0}{\Gamma(k-j)x^{k+j-2}}}}\right)^{1/2}
		< cn.
\end{equation*}
Thus,
\begin{equation} \label{int3}
\textup{E}\left[N\left(1-\frac{\log{\log{n}}}{n}, 1\right)\right] = O(\log{\log{n}}).
\end{equation}

For $(-1, -1+\frac{\log{\log{n}}}{n})$ we will substitute in $-x$.  Notice that since $f(\phi)$ is continuous and non-zero on $[-\pi,\pi]$, we can bound it from below by a constant $\frac{m}{2\pi}>0$.  We then have
\begin{equation*}
A \geq \frac{m}{2\pi}\int_{-\pi}^{\pi}\frac{1-(-x)^{n+1}e^{-i(n+1)\phi}}{1+xe^{-i\phi}} \cdot \frac{1-(-x)^{n+1}e^{i(n+1)\phi}}{1+xe^{i\phi}}d\phi
	=    m\sum_{k=0}^n(-x)^{2k},\\
\end{equation*}
since $f(\phi) \equiv \frac{1}{2\pi}$ in the independent case.  Similarly, we can bound $f(\phi)$ from above by $\frac{M}{2\pi} < \infty$.  This yields the inequality
\begin{flalign*}
C &\leq \frac{M}{2\pi}\int_{-\pi}^{\pi}\frac{-(n+1)(-x)^ne^{-i(n+1)\phi}\left(1 + xe^{-i\phi}\right) + \left(1 - (-x)^{n+1}e^{-i(n+1)\phi}\right)e^{-i\phi}}{\left(1 + xe^{-i\phi}\right)^2}\\
  &  \quad \cdot \frac{-(n+1)(-x)^ne^{i(n+1)\phi}\left(1 + xe^{i\phi}\right) + \left(1 - (-x)^{n+1}e^{i(n+1)\phi}\right)e^{i\phi}}
  {(1 + xe^{i\phi})^2}d\phi\\ 
	&=    M\sum_{k=0}^nk^2(-x)^{2k-2}
	\leq M\frac{n^2}{(-x)^2}\sum_{k=0}^n(-x)^{2k}.
\end{flalign*}
It follows that
\begin{equation*}
\frac{\sqrt{AC-B^2}}{A} < \left(\frac{C}{A}\right)^{1/2}	< cn.
\end{equation*}
Thus,
\begin{equation} \label{int4}
\textup{E}\left[N\left(-1,-1+\frac{\log{\log{n}}}{n}\right)\right] = O(\log{\log{n}}).
\end{equation}
Combining \eqref{int1}-\eqref{int4}, the result then follows.
\end{proof}

\section{Proof of Main Theorem}

Now that we have found an upper bound for the expected number of zeros in our initial four intervals, we will spend the rest of our time deriving explicit values for the intervals $(-1+\frac{\log{\log{n}}}{n}, -1+\elg)$ and $(1-\elg, 1-\frac{\log{\log{n}}}{n})$.  As stated before, this will be done by deriving asymptotic values for the expressions given in \eqref{Kac}.  We will formulate this result as an additional lemma.

\begin{lemma} \label{actual roots}
The expected number of zeros for the polynomial $P_n(x)$ in each of the intervals $(-1+\frac{\log{\log{n}}}{n}, -1+\elg)$ and $(1-\elg, 1-\frac{\log{\log{n}}}{n})$ is given by the following:
\renewcommand{\theenumi}{\roman{enumi}}
\renewcommand{\labelenumi}{(\theenumi)}
\begin{enumerate}
\item For $f\in C\left([-\pi,\pi]\right)$,
\begin{flalign*}
\Exp\left[N\left(-1+\frac{\log{\log{n}}}{n}, -1+\elg\right)\right] &= \Exp\left[N\left(1-\elg, 1-\frac{\log{\log{n}}}{n}\right)\right]\\
	&\sim \frac{1}{2\pi}\log{n}.
\end{flalign*}
\item For $f\in C^1\left([-\pi,\pi]\right)$,
\begin{flalign*}
\Exp\left[N\left(-1+\frac{\log{\log{n}}}{n}, -1+\elg\right)\right] &= \Exp\left[N\left(1-\elg, 1-\frac{\log{\log{n}}}{n}\right)\right]\\
	&= \frac{1}{2\pi}\log{n} + O\left(\log{\log{n}}\right).
\end{flalign*}
\end{enumerate}
\end{lemma}
\begin{proof}
Consider the interval $(1-\elg, 1 - \lgn)$.  For $x=1-y$, define $g(y)=y\frac{\log{n}}{\log{\log{n}}}$.  Also, let $M>0$ be chosen such that $f(\phi)\leq M$, for any $\phi \in [-\pi,\pi]$.  Recalling \eqref{Kac}, for $A$ we have
\begin{equation} \label{a}
\begin{split}
A &= \int_{-\pi}^{\pi}\frac{1}{(1-xe^{-i\phi})(1-xe^{i\phi})}f(\phi)d\phi\\
  &  \quad + x^{n+1}\int_{-\pi}^{\pi} \frac{-\left(e^{-i(n+1)\phi} + e^{i(n+1)\phi}\right) + x^{n+1}}{(1-xe^{-i\phi})(1-xe^{i\phi})}f(\phi)d\phi\\
  &= A_1 + O\left(x^{n+1}A_1\right).\\
\end{split}
\end{equation}
Our next step is to determine the asymptotics for $A_1$.  We have,
\begin{flalign*}  
A_1 &= 2\int_0^{g(y)}\frac{1}{1 - 2x\cos{\phi} + x^2}f(\phi)d\phi + 2\int_{g(y)}^{\pi}\frac{1}{1 - 2x\cos{\phi} + x^2}f(\phi)d\phi.\\
\end{flalign*}
Looking at the first integral yields
\begin{flalign*}
	&= 2\int_{0}^{g(y)}\frac{f(\phi)}{1-2(1-y)\left(1-\frac{\phi^2}{2}+O\left(\phi^4\right)\right)+(1-y)^2}d\phi\\
      &= 2\int_{0}^{g(y)}\frac{f(\phi)}{\phi^2 + y^2}d\phi + O(1).
\end{flalign*}
For $f\in\cont$ this becomes      
\begin{equation*}
			\sim 2\int_{0}^{g(y)}\frac{f(0)}{\phi^2+y^2}d\phi = \frac{2f(0)}{y}\arctan{\left(\frac{g(y)}{y}\right)},\\
\end{equation*}		  
while for $f\in\diff$ we have the more detailed representation
\begin{flalign*}		  
      &= 2\int_{0}^{g(y)}\frac{f(0) + f'(\phi_0)\phi}{\phi^2+y^2}d\phi + O(1) \quad \mbox{(where $\phi_0 \in (0,\phi)$)}\\
		  &= \frac{2f(0)}{y}\arctan{\left(\frac{g(y)}{y}\right)} + O\left(\log{\log{n}}\right).
\end{flalign*}
For the second integral in $A_1$,
\begin{flalign*}
&= \int_{g(y)}^{(g(y))^{1/3}}\frac{2}{1 - 2x\cos{\phi} + x^2}f(\phi)d\phi + \int_{(g(y))^{1/3}}^{\pi}\frac{2}{1 - 2x\cos{\phi} + x^2}f(\phi)d\phi\\
    &\sim \int_{g(y)}^{(g(y))^{1/3}}\frac{2f(0)}{y^2 + \phi^2}d\phi + \int_{(g(y))^{1/3}}^{\pi}\frac{2}{1 - 2x\cos{\phi} + x^2}f(\phi)d\phi\\
    &=   O\left(\frac{1}{g(y)}\right).
\end{flalign*}
It follows that
\begin{flalign*}
A_1 &\sim \frac{2f(0)}{y}\arctan{\left(\frac{g(y)}{y}\right)}, \qquad \mbox{for $f\in\cont$},\\
A_1 &= \frac{2f(0)}{y}\arctan{\left(\frac{g(y)}{y}\right)} + O\left(\frac{1}{g(y)}\right), \qquad \mbox{for $f\in\diff$}.
\end{flalign*}
Also, notice that 
\begin{equation*}
x^{n+1}A_1 \leq \left(1-\lgn\right)^{n+1}\frac{c}{y} = o\left(\frac{1}{g(y)}\right).
\end{equation*}
Combining these results and plugging into \eqref{a} yields
\begin{equation} \label{A}
\begin{split}
A &\sim \frac{2f(0)}{y}\arctan{\left(\frac{g(y)}{y}\right)}, \qquad \mbox{for $f\in\cont$},\\
A &= \frac{2f(0)}{y}\arctan{\left(\frac{g(y)}{y}\right)} + O\left(\frac{1}{g(y)}\right), \qquad \mbox{for $f\in\diff$}.
\end{split}
\end{equation}

Considering $B$ next, 
\begin{equation} \label{b}
\begin{split}
B &= \int_{-\pi}^{\pi}\frac{e^{i\phi}}{(1-xe^{-i\phi})(1-xe^{i\phi})^2}f(\phi)d\phi\\
	&  \quad + \int_{-\pi}^{\pi}\frac{-(n+1)x^ne^{i(n+1)\phi}\left(1 - xe^{i\phi}\right)}{(1-xe^{-i\phi})(1-xe^{i\phi})^2}f(\phi)d\phi\\
	&  \quad + x^{n+1}\int_{-\pi}^{\pi} \frac{(n+1)x^n\left(1 - xe^{i\phi}\right) - e^{-in\phi}\left(1 - x^{n+1}e^{i(n+1)\phi}\right) - e^{i(n+2)\phi}}{(1-xe^{-i\phi})(1-xe^{i\phi})^2}f(\phi)d\phi\\
	&= B_1 + B_2 + O\left(x^{n+1}\left(|B_1| + |B_2|\right)\right).
\end{split}	
\end{equation}
Since we will end up showing that the $B_1$ term dominates, we can rewrite this as
\begin{equation*} 
B = B_1 + B_2 + O\left(x^{n+1}B_1\right).
\end{equation*}
To analyze $B_1$ we will split it into two integrals,
\begin{flalign*}
B_1 &= 2\int_{0}^{g(y)}\frac{\cos{\phi} - 1 +y}{\left(1 - 2(1 - y)\cos{\phi} + (1-y)^2\right)^2}f(\phi)d\phi\\
		&  \quad + 2\int_{g(y)}^{\pi}\frac{\cos{\phi} - 1 + y}{\left(1 - 2(1 - y)\cos{\phi} + (1-y)^2\right)^2}f(\phi)d\phi.\\
\end{flalign*}
For the first integral we have
\begin{flalign*}		
 		&= 2\int_{0}^{g(y)}\frac{y-\frac{\phi^2}{2}+O\left(\phi^4\right)}{\left(1-2(1-y)\left(1-\frac{\phi^2}{2}+O\left(\phi^4\right)\right) + (1-y)^2\right)^2}f(\phi)d\phi\\
    &= 2\int_{0}^{g(y)}\frac{y}{\left(\phi^2+y^2\right)^2}f(\phi)d\phi + O\left(\frac{1}{y}\right).\\
\end{flalign*}
When $f\in\cont$ this becomes
\begin{equation*}    
    \sim 2\int_{0}^{g(y)}\frac{yf(0)}{\left(\phi^2+y^2\right)^2}d\phi \sim \frac{f(0)}{y^2}\arctan{\left(\frac{g(y)}{y}\right)},
\end{equation*}
while for $f\in\diff$
\begin{flalign*}
		&= 2\int_{0}^{g(y)}\frac{yf(0) + yf'(\phi_0)\phi}{\left(\phi^2+y^2\right)^2}d\phi + O\left(\frac{1}{y}\right)
		  \quad \left(\mbox{where $\phi_0\in(0,\phi)$}\right)\\
    &= \frac{f(0)}{y^2}\arctan{\left(\frac{g(y)}{y}\right)} + O\left(\frac{1}{yg(y)}\right).
\end{flalign*}
For the second  integral in $B_1$,
\begin{flalign*}
&        \left|2\int_{g(y)}^{\pi}\frac{\cos{\phi} - x}{(1 - 2x\cos{\phi} + x^2)^2}f(\phi)d\phi\right|\\
&\leq    2\int_{g(y)}^{(g(y))^{1/3}}\frac{\left|\cos{\phi} - x\right|f(\phi)}{(1 - 2x\cos{\phi} + x^2)^2}d\phi + 																 							 2\int_{(g(y))^{1/3}}^{\pi}\frac{\left|\cos{\phi} - x\right|f(\phi)}{(1 - 2x\cos{\phi} + x^2)^2}d\phi\\
&\sim 2f(0)\int_{g(y)}^{(g(y))^{1/3}}\frac{\left|y - \frac{\phi^2}{2} + O\left(\phi^4\right)\right|}{\left(\phi^2+y^2\right)^2}d\phi + 																 2\int_{(g(y))^{1/3}}^{\pi}\frac{\left|\cos{\phi} - x\right|f(\phi)}{(1 - 2x\cos{\phi} + x^2)^2}d\phi\\
&=       O\left(\frac{y}{(g(y))^3}\right).
\end{flalign*}
Thus,
\begin{equation*}
\begin{split}
B_1 &\sim \frac{f(0)}{y^2}\arctan{\left(\frac{g(y)}{y}\right)}, \qquad \mbox{for $f\in\cont$},\\
B_1 &= \frac{f(0)}{y^2}\arctan{\left(\frac{g(y)}{y}\right)} + O\left(\frac{1}{yg(y)}\right), \qquad \mbox{for $f\in\diff$}.
\end{split}
\end{equation*}

Looking at $B_2$ next,
\begin{flalign*}
|B_2| &= \left|\int_{-\pi}^{\pi}\frac{-(n+1)x^ne^{i(n+1)\phi}(1-xe^{i\phi})}{(1-xe^{-i\phi})(1-xe^{i\phi})^2}f(\phi)d\phi\right|\\
		&\leq \int_{-\pi}^{\pi}\frac{(n+1)x^n}{(1-xe^{-i\phi})(1-xe^{i\phi})}f(\phi)d\phi\\
		&\sim c(n+1)\frac{(1-y)^n}{y},\\
\end{flalign*}
where the last line follows from our work on $A$.  Now, notice that for $y=\lgn$ we have
\[
(n+1)\frac{(1-y)^n}{y} \sim \frac{n^2}{\log{n}\log{\log{n}}} = \frac{1}{yg(y)}.
\]
Since $(1-y)^ny$ is a decreasing function as $y$ increases, it follows that for $y\in(\lgn,\elg)$
\begin{equation} \label{ny}
(n+1)(1-y)^ny \leq \frac{y}{g(y)} \Rightarrow (n+1)\frac{(1-y)^n}{y} = O\left(\frac{1}{yg(y)}\right).
\end{equation}
This also implies that
\begin{equation*}
x^{n+1}(B_1 + B_2) < c(n+1)\frac{(1-y)^{n+1}}{y}.
\end{equation*}
Thus,
\begin{equation} \label{B}
\begin{split}
B &\sim \frac{f(0)}{y^2}\arctan{\left(\frac{g(y)}{y}\right)}, \qquad \mbox{for $f\in\cont$},\\
B &= \frac{f(0)}{y^2}\arctan{\left(\frac{g(y)}{y}\right)} + O\left(\frac{1}{yg(y)}\right), \qquad \mbox{for $f\in\diff$}.
\end{split}
\end{equation}

Turning now to $C$, 
\begin{equation} \label{c}
\begin{split}
C &= \int_{-\pi}^{\pi}\frac{1}{(1-xe^{-i\phi})^2(1-xe^{i\phi})^2}f(\phi)d\phi + \int_{-\pi}^{\pi}\frac{(n+1)^2x^{2n}}{(1 - xe^{-i\phi})(1 - xe^{i\phi})}f(\phi)d\phi\\
	&  \quad + 2\int_{-\pi}^{\pi}\frac{-(n+1)x^ne^{-in\phi}\left(1-x^{n+1}e^{i(n+1)\phi}\right)}{(1-xe^{-i\phi})(1-xe^{i\phi})^2}f(\phi)d\phi\\
	&  \quad + x^{n+1}\int_{-\pi}^{\pi}\frac{-\left(e^{-i(n+1)\phi} + e^{i(n+1)\phi}\right) + x^{n+1}}{(1-xe^{-i\phi})^2(1-xe^{i\phi})^2}f(\phi)d\phi\\
	&= C_1 + C_2 + C_3 + O\left(x^{n+1}C_1\right).
\end{split}
\end{equation}
Splitting up $C_1$ gives us
\begin{equation*} 
C_1 =			 2\int_{0}^{g(y)}\frac{1}{(1 - 2x\cos{\phi} + x^2)^2}f(\phi)d\phi + 2\int_{g(y)}^{\pi}\frac{1}{(1 - 2x\cos{\phi} + x^2)^2}f(\phi)d\phi.
\end{equation*}
For the first term we have
\begin{flalign*}
		 &= 2\int_0^{g(y)} \frac{1}{\left(1 - 2(1 - y)\left(1 - \frac{\phi^2}{2} + O\left(\phi^4\right)\right) + (1- y)^2\right)^2}f(\phi)d\phi\\
		 &= 2\int_0^{g(y)} \frac{f(\phi)}{\left(\phi^2 + y^2\right)^2}d\phi + O\left(\frac{1}{y^2}\right).\\
\end{flalign*}		 
When $f\in\cont$ it follows that
\begin{equation*}
		 \sim 2\int_0^{g(y)} \frac{f(0)}{\left(\phi^2 + y^2\right)^2}d\phi \sim \frac{f(0)}{y^3}\arctan{\left(\frac{g(y)}{y}\right)},\\
\end{equation*}
while for $f\in\diff$
\begin{flalign*}
&= 2\int_0^{g(y)} \frac{f(0) + f'(\phi_0)\phi}{\left(\phi^2 + y^2\right)^2}d\phi + O\left(\frac{1}{y^2}\right) \quad \mbox{(where $\phi_0 \in(0,\phi)$)}\\
     &= \frac{f(0)}{y^3}\arctan{\left(\frac{g(y)}{y}\right)} + O\left(\frac{1}{y^2g(y)}\right).\\
\end{flalign*}
For the second integral in $C_1$ we have
\begin{flalign*}
    &=       2\int_{g(y)}^{(g(y))^{1/3}}\frac{f(\phi)}{\left(1 - 2x\cos{\phi} + x^2\right)^2}d\phi + 2\int_{(g(y))^{1/3}}^{\pi}\frac{f(\phi)}{\left(1 - 2x\cos{\phi} + x^2\right)^2}d\phi\\
    &\sim 2f(0)\int_{g(y)}^{(g(y))^{1/3}}\frac{1}{\left(\phi^2+y^2\right)^2}d\phi + 2\int_{(g(y))^{1/3}}^{\pi}\frac{f(\phi)}{\left(1 - 2x\cos{\phi} + x^2\right)^2}d\phi\\
    &=       O\left(\frac{1}{(g(y))^3}\right).
\end{flalign*}
Thus,
\begin{equation} \label{C1}
\begin{split}
C_1 &\sim \frac{f(0)}{y^3}\arctan{\left(\frac{g(y)}{y}\right)}, \qquad \mbox{for $f\in\cont$},\\
C_1 &= \frac{f(0)}{y^3}\arctan{\left(\frac{g(y)}{y}\right)} + O\left(\frac{1}{y^2g(y)}\right), \qquad \mbox{for $f\in\diff$}.
\end{split}
\end{equation}

For $C_2$,
\begin{equation*}
C_2 =       \int_{-\pi}^{\pi}\frac{(n+1)^2x^{2n}}{(1-xe^{-i\phi})(1-xe^{i\phi})}f(\phi)d\phi
    \sim c(n + 1)^2\frac{(1 - y)^{2n}}{y},\\
\end{equation*}
while for $C_3$ we have
\begin{flalign*}
|C_3| &\leq c\int_{-\pi}^{\pi}\left|\frac{-(n+1)x^ne^{-in\phi}}{(1-xe^{-i\phi})(1-xe^{i\phi})^2}f(\phi)\right|d\phi\\
      &\leq c(C_1C_2)^{1/2} \quad \mbox{(by Cauchy-Schwarz)}\\
      &\sim c\frac{(n+1)(1-y)^n}{y^2}.
\end{flalign*}
Also, note that
\[
x^{n+1}C_1 \sim c\frac{(1-y)^{n+1}}{y^3} \leq \frac{(n+1)(1-y)^n}{y^2}.
\]
From \eqref{ny} we know that $(n+1)(1-y)^n = O\left(1/g(y)\right)$, which implies
\begin{equation} \label{remainder}
\begin{split}
\frac{(n+1)(1-y)^n}{y^2} &= O\left(\frac{1}{y^2g(y)}\right),\\
\frac{(n+1)^2(1-y)^{2n}}{y} &= O\left(\frac{1}{y g^2(y)}\right).
\end{split}
\end{equation}
Thus, we can conclude that
\begin{equation} \label{C}
\begin{split}
C &\sim \frac{f(0)}{y^3}\arctan{\left(\frac{g(y)}{y}\right)}, \qquad \mbox{for $f\in\cont$},\\
C &= \frac{f(0)}{y^3}\arctan{\left(\frac{g(y)}{y}\right)} + O\left(\frac{1}{y^2g(y)}\right), \qquad \mbox{for $f\in\diff$}.
\end{split}
\end{equation}

Combining \eqref{A}, \eqref{B}, and \eqref{C}, for $f\in\cont$ we then have the expression
\begin{equation} \label{abccont}
\frac{\sqrt{AC - B^2}}{A} \sim \frac{1}{2y},
\end{equation}
while for $f\in\diff$ this becomes
\begin{equation} \label{abc}
\frac{\sqrt{AC - B^2}}{A} = \frac{\sqrt{\frac{f^2(0)}{y^4}\left[\arctan{\left(\frac{g(y)}{y}\right)}\right]^2 + O\left(\frac{1}{y^3g(y)}\right)}}
			 {\frac{2f(0)}{y}\arctan{\left(\frac{g(y)}{y}\right)} + O\left(\frac{1}{g(y)}\right)}
			= \frac{1}{2y} + O\left(\frac{1}{g(y)}\right).
\end{equation}
Plugging these into \eqref{Kac} gives us 
\begin{equation} \label{positive roots cont}
\textup{E}\left[N\left(1-\elg, 1-\lgn\right)\right] \sim \frac{1}{2\pi}\log{n},
\end{equation}	
for $f\in\cont$, and
\begin{equation} \label{positive roots diff}
\textup{E}\left[N\left(1-\elg, 1-\lgn\right)\right] = \frac{1}{2\pi}\log{n} + O\left(\log{\log{n}}\right),
\end{equation}	
for $f\in\diff$.

To handle the interval from $(-1+\lgn, -1+\elg)$ we will substitute in $-x=-1+y$, where $x\in (1-\elg,1-\lgn)$.  We then have
\begin{equation} \label{ABC-}
\begin{split}
A &=       \int_{-\pi}^{\pi}\frac{1-(-x)^{n+1}e^{-i(n+1)\phi}}{1+xe^{-i\phi}} \cdot \frac{1-(-x)^{n+1}e^{i(n+1)\phi}}{1+xe^{i\phi}}f(\phi)d\phi,\\
B &= \int_{-\pi}^{\pi}\left(\frac{1-(-x)^{n+1}e^{-i(n+1)\phi}}{1+xe^{-i\phi}}\right)\\
  &  \quad \cdot \left(\frac{-(n+1)(-x)^ne^{i(n+1)\phi}\left(1 + xe^{i\phi}\right) + \left(1 - (-x)^{n+1}e^{i(n+1)\phi}\right)e^{i\phi}}{\left(1 + xe^{i\phi}\right)^2}\right)f(\phi)d\phi,\\
C &= \int_{-\pi}^{\pi}\left(\frac{-(n+1)(-x)^ne^{-i(n+1)\phi}\left(1 + xe^{-i\phi}\right) + \left(1 - (-x)^{n+1}e^{-i(n+1)\phi}\right)e^{-i\phi}} {\left(1 + xe^{-i\phi}\right)^2}\right)\\
  &  \quad \cdot \left(\frac{-(n+1)(-x)^ne^{i(n+1)\phi}\left(1 + xe^{i\phi}\right) + \left(1 - (-x)^{n+1}e^{i(n+1)\phi}\right)e^{i\phi}} {\left(1 + xe^{i\phi}\right)^2}\right)f(\phi)d\phi.\\
\end{split}  
\end{equation}
Using the fact that $f(\pi)=f(-\pi)$, we can apply methods similar to those from the positive case to derive the formulas
\begin{equation} \label{ABCcont-}
\begin{split}
A &\sim \frac{2f(\pi)}{y}\arctan{\left(\frac{g(y)}{y}\right)},\\
B &\sim -\frac{f(\pi)}{y^2}\arctan{\left(\frac{g(y)}{y}\right)},\\
C &\sim \frac{f(\pi)}{y^3}\arctan{\left(\frac{g(y)}{y}\right)},\\
\end{split}
\end{equation}
for $f\in\cont$, and
\begin{equation} \label{ABCdiff-}
\begin{split}
A &= \frac{2f(\pi)}{y}\arctan{\left(\frac{g(y)}{y}\right)} + O\left(\frac{1}{g(y)}\right),\\
B &= -\frac{f(\pi)}{y^2}\arctan{\left(\frac{g(y)}{y}\right)} + O\left(\frac{1}{yg(y)}\right),\\
C &= \frac{f(\pi)}{y^3}\arctan{\left(\frac{g(y)}{y}\right)} + O\left(\frac{1}{y^2g(y)}\right),
\end{split}
\end{equation}
for $f\in\diff$.
Combining \eqref{ABCcont-} and \eqref{ABCdiff-}, we then have the same expressions as before.  That is, \eqref{ABCcont-} gives
\begin{equation} \label{abc-cont}
\frac{\sqrt{AC - B^2}}{A} = \frac{1}{2y},
\end{equation}
while from \eqref{ABCdiff-} we again have
\begin{equation} \label{abc-}
\frac{\sqrt{AC - B^2}}{A} = \frac{\sqrt{\frac{f^2(\pi)}{y^4}\left[\arctan{\left(\frac{g(y)}{y}\right)}\right]^2 + O\left(\frac{1}{y^3g(y)}\right)}}
			{\frac{2f(\pi)}{y}\arctan{\left(\frac{g(y)}{y}\right)} + O\left(\frac{1}{g(y)}\right)}
			= \frac{1}{2y} + O\left(\frac{1}{g(y)}\right).
\end{equation}
Plugging these into \eqref{Kac} gives us 
\begin{equation} \label{negative roots cont}
\textup{E}\left[N\left(-1+\lgn, -1+\elg\right)\right] \sim \frac{1}{2\pi}\log{n},
\end{equation}
for $f\in\cont$, and
\begin{equation} \label{negative roots diff}
\textup{E}\left[N\left(-1+\lgn, -1+\elg\right)\right] = \frac{1}{2\pi}\log{n} + O\left(\log{\log{n}}\right),
\end{equation}
for $f\in\diff$.

\end{proof}
We are now able to prove Theorem \ref{roots}.
\begin{proof}[Proof of Theorem \ref{roots}]
Combining the results of Lemmas \ref{bound} and \ref{actual roots}, we have
\begin{equation} \label{-1to1}
\begin{split}
\textup{E}[N(-1, 1)] &\sim \frac{1}{\pi}\log{n}, \qquad \mbox{for $f\in\cont$},\\
\textup{E}[N(-1, 1)] &=    \frac{1}{\pi}\log{n} + O\left(\log{\log{n}}\right), \qquad \mbox{for $f\in\cont$}.\\
\end{split}
\end{equation}
Recalling from Section 1 that $\textup{E}[N(-\infty,\infty)] = 2\textup{E}[N(-1, 1)]$, the result then follows.
\end{proof}

\section{Acknowledgments}
The author would like to thank his thesis advisor, Professor Michael Cranston, for countless helpful discussions with the subject.  His advice and suggestions were invaluable.  The author would also like to acknowledge Professor Stanislav Molchanov for suggesting the problem and the idea of using the spectral density of the covariance function, as well as the referee for his insightful comments.  Finally, a thank you is owed to Mr. Phillip McRae for reading through a copy of this manuscript.

\end{document}